\documentclass{article}
\usepackage{hyperref}
\hypersetup{breaklinks=true}
\usepackage{relsize}
\usepackage{graphicx}
\usepackage{amsmath, amsthm, amssymb, amscd, mathrsfs, eucal, epsfig,color}
\usepackage[english]{babel}

\usepackage{pinlabel}
\usepackage{pgf}

\newtheorem{theorem}{Theorem}[section]
\newtheorem{corollary}[theorem]{Corollary}
\newtheorem{lemma}[theorem]{Lemma}

\newtheorem*{prop*}{Proposition}
\newtheorem{example}[theorem]{Example}
\newtheorem*{theorem*}{Theorem}
\newtheorem*{corollary*}{Corollary}
\newtheorem*{lemma*}{Lemma}

\theoremstyle{definition}

\theoremstyle{remark}

\newlength{\larg}
\newcommand{\Mod}{\textrm{Mod}}

\begin{document}
\title{Big mapping class groups are not acylindrically hyperbolic}
\author{Juliette Bavard\footnote{Supported by grants from R\'{e}gion Ile-de-France.}, Anthony Genevois}
\date{\today}

\maketitle
\begin{abstract} 
We give a criterion to prove that some groups are not acylindrically hyperbolic. As an application, we prove that the mapping class group of an infinite type surface is not acylindrically hyperbolic.
\end{abstract}

\setcounter{tocdepth}{2}
%\tableofcontents

\section*{Introduction}

%\subsection*{Actions of mapping class groups on Gromov-hyperbolic spaces}
If $S$ is an orientable surface, the \emph{mapping class group of $S$}, denoted by $\Mod(S)$, is the group of orientation-preserving homeomorphisms of $S$ up to isotopy. If $S$ has finite genus and contains only finitely many punctures and boundary components, we say that $S$ is of \emph{finite type}. If not, we say that $S$ is an \emph{infinite type surface}. 

%\subsubsection*{Finite type surfaces}
Mapping class groups of finite type surfaces has been studied for many years. If $S$ has enough complexity (i.e., $S$ is not a torus with less than two punctures nor a sphere with less than five punctures), then $\Mod (S)$ acts by isometries on its \emph{curve complex} $\mathcal{C}(S)$. The curve complex of $S$ is a simplicial space whose vertices are isotopy classes of essential simple closed curves in $S$, and whose simplexes of dimension $n$ are the $n$-uplets of vertices which have disjoint representatives in $S$. This complex is known to be Gromov-hyperbolic by a theorem of Howard Masur and Yair Minsky (see \cite{Masur-Minsky}). Moreover, Mladen Bestvina and Koji Fujiwara have proved in \cite{Bestvina-Fujiwara} that the action of $\Mod(S)$ on $\mathcal{C}(S)$ is \emph{weakly properly discontinuous}, and Bowditch has proved more generally that this action is \emph{acylindrical} in \cite{Bowditch}. In the langage introduced by Osin in \cite{Osin1}, it follows that $\Mod(S)$ is \emph{acylindrically hyperbolic} (see Section \ref{1} for a precise definition). In particular, this implies that: 
\begin{itemize}
	\item $\Mod(S)$ is \emph{SQ-universal}, i.e., any countable group embeds in a quotient of $\Mod(S)$. In particular, $\Mod(S)$ contains non abelian free subgroups and it contains uncountably many normal subgroups (see \cite{DGO});
	\item the dimension of the second bounded cohomology of $\Mod(S)$ is infinite (see \cite{Bestvina-Fujiwara});
	\item any asymptotic cone of $\Mod(S)$ contains a cut point (\cite{Sisto1}, see also \cite{Behrstock}).
\end{itemize}
More specifically, this approach leads Fran\c{c}ois Dahmani, Vincent Guirardel and Denis Osin to answer several open questions on mapping class groups in \cite{DGO}. Motivating by these results, it is natural to ask whether the same ideas apply to infinite type surfaces.

%\subsubsection*{Infinite type surfaces}
Curve complexes of infinite type surfaces have finite diameter so that they are not useful to study the associated mapping class groups from the view point of coarse geometry. For the specific case of the plane minus a Cantor set, it has been recently proved by the first author in \cite{Ju} that the mapping class group acts by isometries on a Gromov-hyperbolic graph of infinite diameter. This graph is the \emph{ray graph}, defined by Danny Calegari in \cite{blog-Calegari}. Even more recently, Javier Aramayona, Ariadna Fossas and Hugo Parlier have extended this result to more infinite type surfaces (see \cite{Aramayona-and-co}). However, these actions are far from being acylindrical or from containing WPD isometries: for every finite set $V$ of vertices of the graph, there are infinitely many elements of the mapping class group which fix $V$. Thus, although the ray graph seems to be the natural analogous of the curve graph in the context of infinite type surfaces, a natural question would be to determine if there exists a more suitable hyperbolic graph on which the associated mapping class group acts acylindrically. The main result of this article answers negatively this question. Namely, for any infinite type surface $S$, we prove that every acylindrical action of $\Mod(S)$ on a hyperbolic space is necessary elementary. Equivalently: 

\begin{theorem*}\textbf{\textrm{\ref{theo:big_Mod}}}
If $S$ is an infinite type surface, then the group $\Mod(S)$ is not acylindrically hyperbolic.
\end{theorem*}

\subsection*{Acknowledgment} The first author would like to thank her PhD advisor Fr\'{e}d\'{e}ric Le Roux for many advice and explanations, especially about the classification of infinite type surfaces and the space of ends.

\section{A criterion of non acylindrical hyperbolicity}\label{1}

\noindent
The action by isometries of a group $G$ on a metric space $X$ is \emph{acylindrical} provided that, for every $d \geq 0$, there exist $R,N \geq 0$ such that for any $x,y \in X$, 
\begin{center}
$d(x,y) \geq R \Rightarrow \# \{ g \in G \mid d(x,gx),d(y,gy) \leq d \} \leq N$.
\end{center}
Furthermore, if $X$ is hyperbolic, we say that the action is \emph{non elementary} if the limit set $\Lambda(G)$ of $G$ is infinite. A group is \emph{acylindrically hyperbolic} if it admits a non elementary acylindrical action on a hyperbolic space.

\medskip \noindent
The main result of this section is:

\begin{theorem}\label{critère}
Let $G$ be a group with a collection of subgroups $\mathcal{H}$ satisfying:
\begin{itemize}
	\item[(i)] for every $H \in \mathcal{H}$ and $h \in H$, the centralizer of $h$ is not virtually cyclic,
	\item[(ii)] for every $g \in G$ and $N \geq 1$, there exist $n \geq N$ and $H \in \mathcal{H}$ such that $H \cap H^{g^{n}}$ is infinite.
\end{itemize}
Then $G$ is not acylindrically hyperbolic.
\end{theorem}

\noindent
Therefore, given a group $G$ satisfying the hypotheses of our theorem, we want to prove that any acylindrical action of $G$ on a hyperbolic space turns out to be elementary. For this purpose, the following trichotomy will be useful (see \cite[Theorem 1.1]{Osin1}):

\begin{theorem}\label{Osin1}
Let $G$ be a group acting acylindrically on a hyperbolic space $X$. Then exactly one of the following situation happens:
\begin{itemize}
	\item[(i)] $G$ contains infinitely many pairwise independent loxodromic isometries,
	\item[(ii)] $G$ is virtually cyclic and contains a loxodromic isometry,
	\item[(iii)] the action $G \curvearrowright X$ has a bounded orbit.
\end{itemize} 
In particular, the action is non elementary if and only if (i) holds.
\end{theorem}

\noindent
Recall that two isometries of a hyperbolic space are \textit{independent} if their points at infinity are pairwise distinct. As a consequence, by applying the previous classification to the cyclic subgroups of $G$, we deduce that any element $g \in G$ is either loxodromic or elliptic:

\begin{corollary}\label{Osin2}
If a group acts acylindrically on a hyperbolic space, none of its elements are parabolic. 
\end{corollary}

\noindent
A last result needed to prove our theorem is \cite[Corollary 6.6]{DGO}:

\begin{lemma}\label{DGO}
Let $G$ be a group acting acylindrically on a hyperbolic space $X$. If $g \in G$ induces a loxodromic isometry of $X$ then its centralizer is virtually cyclic.
\end{lemma}

\noindent
\textbf{Proof of Theorem \ref{critère}.} Given a group $G$ satisfying the hypotheses of our theorem, we want to prove that any acylindrical action of $G$ on a $\delta$-hyperbolic space $X$ is elementary. First, it clearly follows from Lemma \ref{DGO} that any element of a subgroup of $\mathcal{H}$ does not induce a loxodromic isometry of $X$: according to Corollary \ref{Osin2}, it induces an elliptic isometry.

\medskip \noindent
Suppose by contradiction that $G$ contains a loxodromic isometry $g$, and fix a constant $A >0$. Because $g$ is loxodromic, there exists some $n \geq 1$ such that $n\|g\| = \|g^n \| \geq A$, where $\| \cdot \|$ denotes the translation length \cite[Proposition 10.6.1]{CDP}. By our hypotheses, there exist some $H \in \mathcal{H}$ and $k \geq n$ such that $H \cap H^{g^k}$ is infinite. Notice that $k \geq n$ implies $\| g^k \| \geq A$ as well. As we have seen, the induced action $H \curvearrowright X$ does not contain any loxodromic isometry, so it follows from Theorem \ref{Osin1} that the orbits of $H$ are bounded; therefore, there exists a point $x_0 \in X$ whose orbit under $H$ has diameter at most $5 \delta$ \cite[Lemma III.$\Gamma$.3.3]{BH}. Notice that, for all $h \in H \cap H^{g^k}$, we have $d(x_0,hx_0) \leq 5 \delta$ and
\begin{center}
$d(g^kx_0,hg^kx_0)=d(x_0,g^{-k}hg^k x_0) \leq 5 \delta$
\end{center}
since $g^{-k}hg^k \in g^{-k}H^{g^k}g^k=H$. Therefore, we have proved that, for any $A>0$, there exist $x,y \in X$ satisfying $d(x,y) \geq A$ and such that the set $\{ g \in G \mid d(x,gx),d(y,gy) \leq 5 \delta \}$ is infinite. This contradicts the acylindricity of the action $G \curvearrowright X$. Thus, this action does not contain any loxodromic isometry, and we deduce that it is necessarily elementary. $\square$

\begin{example}
As a first application, we claim that $SO(n)$ is not acylindrically hyperbolic for every $n \in \mathbb{N}$. Indeed, let $\mathcal{H}$ denote the collection of the subgroups of $SO(n)$ which have the form:
\begin{center}
$PE_kP^{-1} := P \left( \begin{matrix} \mathrm{Id}_k & & \\ & SO(2) & \\ & & \mathrm{Id}_{n-k-2} \end{matrix} \right) P^{-1}$,
\end{center}
where $0 \leq k \leq n$ and $P \in \mathrm{GL}(n)$. Notice that any element of $PE_kP^{-1}$ commutes with any element of $PE_{h}P^{-1}$ provided that $|k-h| \geq 2$; in particular, the centralizer of any element of $PE_kP^{-1}$ is not virtually cyclic if $n \geq 4$. Then, notice that a matrix $M \in SO(n)$ can be written as $P \left( \begin{matrix} A_1 & & 0 \\ & \ddots & \\ 0 & & A_n \end{matrix} \right) P^{-1}$ for some $A_i \in SO(2)$ or $\{ \pm 1 \}$. Up to a conjugacy, we may suppose without loss of generality that $M = \left( \begin{matrix} A_1 & & 0 \\ & \ddots & \\ 0 & & A_n \end{matrix} \right)$. We have  
\begin{center}
$M \left( \begin{matrix} SO(2) &  & & 0 \\ \\ & 1 & & \\ & & \ddots & \\ 0 & & & 1 \end{matrix} \right) M^{-1} = \left( \begin{matrix} SO(2) & &  & 0 \\ \\ & 1 & & \\ & & \ddots & \\ 0 & & & 1 \end{matrix} \right),$
\end{center}
that is to say $E_0^M = E_0$. We conclude that Theorem \ref{critère} applies so that $SO(n)$ is not acylindrically hyperbolic for any $n \geq 4$. Finally, notice that $SO(2)$ and $SO(3)$ are not acylindrically hyperbolic as well since they are respectively abelian and simple (see \cite[Theorem 2.29]{DGO}), so $SO(n)$ is not acylindrically hyperbolic for any $n \in \mathbb{N}$. 
\end{example}

\section{Application to big mapping class groups}

As an application of Theorem \ref{critère}, we get the following theorem, where $\Mod(S)$ denotes the group of orientation-preserving homeomorphisms of $S$ up to isotopy.

\begin{theorem}\label{theo:big_Mod}
Let $S$ be an infinite type surface. Then $\Mod(S)$ is not acylindrically hyperbolic.
\end{theorem}

Before the proof we recall necessary elements on the classification of infinite type surfaces, which can be found in \cite{Richards}; see also \cite{Ghys}. 
\paragraph{Space of ends.}
If $S$ is a non compact surface, consider a sequence $(S_i)$ of compact subsurfaces of $S$ such that:
\begin{itemize}
\item for every $i$, $S_i$ is included into $S_{i+1}$;
\item $\bigcup_i S_i =S$
\end{itemize}
Then every sequence $(U_i)$ satisfying:
\begin{itemize}
\item for every $i$, $U_i$ is a connected component of $S\setminus S_i$;
\item for every $i$, $U_{i+1}$ is included in $U_{i}$;
\end{itemize}
is called an \emph{end} of $S$. It is worth noticing that there are several different topological types of ends. One end $(U_i)$ is said to be:
\begin{enumerate}
\item \emph{isolated} if all but finitely many $U_i$'s are homeomorphic to a cylinder;
\item \emph{accumulated by ends} if every $U_i$ belongs to another end (seen as a sequence);
\item \emph{accumulated by genus} if every $U_i$ has infinite genus.
\end{enumerate}

The \emph{space of ends} (or \emph{ideal boundary}) of the surface $S$ will be denoted by $K$. It is equipped with a natural topology for which it is totally disconnected, separable and compact. Furthermore, this space does not depend on our choice of the sequence $(S_i)$: there exists a natural homeomorphism between the spaces of ends associated to any two such sequences.

\paragraph{Notations.} For convenience, the set of ends accumulated by ends will be denoted by $K_1$, and the set of ends accumulated by genus by $K_2$. Note that both are compact subsets of $K$, and that if $S$ is of infinite type, at least one of them is non empty. These subspaces are of particular interest since, according to a theorem of Ker\'ekj\'art\'o (see \cite{Richards} or \cite{Ghys}), an orientable surface of infinite type is uniquely determined, up to homeomorphism, by its genus and the triple of topological spaces $(K,K_1,K_2)$. If $S'$ is a subsurface of $S$, we say that an end \emph{$(U_i) \in K$ belongs to $S'$} or that \emph{$S'$ contains $(U_i)$} if $U_i \subset S'$ for all but finitely many $i$'s.

\paragraph{Back to the proof of the theorem.} We now use the previous definitions and notations to prove Theorem \ref{theo:big_Mod}.
\begin{proof} If $\Omega$ is a subsurface of $S$, we say that $g\in \Mod(S)$ is \emph{supported in $\Omega$} if $g$ has a representative which is supported in $\Omega$. With the notations we have just defined for the space of ends and its subsets, we consider the two following cases.\\

 \textbf{$(1)$ First case: if $S$ is a surface of finite genus with infinitely many ends ($K_1\neq \emptyset$).} Consider the set $X_1$ of the topological punctured closed disks $\Omega$ of $S$ (up to isotopy) containing at least one end of $K_1$ whose complement in $S$ contains at least two ends of $K$ or has positive genus. For every $\Omega \in X_1$, we denote by $H_\Omega$ the subgroup of $\Mod(S)$ consisting in the homeomorphisms supported in $\Omega$. Finally, we denote by $\mathcal{H}_1$ the collection $\{H_\Omega\}_{\Omega \in X_1}$.\\

  \textbf{$(2)$ Second case: if $S$ is a surface of infinite genus ($K_2\neq \emptyset$).} Consider the set $X_2$ of the subsurfaces $\Omega$ of $S$ (up to isotopy) containing at least one end in $K_2$ whose complement in $S$ contains at least two ends of $K$ or has positive genus. Again, for every such $\Omega$, we denote by $H_\Omega$ the subgroup of $\Mod(S)$ consisting in the homeomorphisms supported in $\Omega$. Finally, we denote by $\mathcal{H}_2$ the collection $\{H_\Omega\}_{\Omega \in X_2}$.\\

\noindent 
Choose $i \in \{1,2\}$ such that $K_i$ (and hence $X_i$) is non empty. We have:
\begin{itemize}
\item For every $\Omega \in X_i$ and $h \in H_\Omega$, every element of $\Mod(S)$ which is supported in $S-\Omega$ commutes with $h$: it follows that the centralizer of $h$ is not virtually cyclic.
\item Let $g \in \Mod(S)$. Note that $\langle g \rangle$ acts on $K_i$. Let $x \in K_i$ be an accumulation point for this action; such an $x$ exists, because $K_i$ is compact. Choose an $\Omega \in X_i$ containing $x$. For every $N\geq 1$, there exists $n\geq N$ such that $g^n(\Omega)\cap \Omega$ contains $x$. 
\begin{enumerate}
\item Case $i=1$: since $x$ is accumulated by ends, there are infinitely many points of $K$ in $g^n(\Omega)\cap \Omega$.
\item Case $i=2$: since $x$ is accumulated by genus, there are infinitely many elements of $\Mod(S)$ which are supported in $g^n(\Omega)\cap \Omega$.
\end{enumerate}
In both cases, $H_\Omega \cap H_\Omega^{g^n}$ is equal to $H_{g^n(\Omega)\cap \Omega}$,  which is infinite.
\end{itemize}
It follows that the collection $\mathcal{H}_i$ satisfies the hypotheses of Theorem \ref{critère}: therefore, $\Mod(S)$ is not acylindrically hyperbolic.
\end{proof}

\addcontentsline{toc}{section}{References}

\bibliographystyle{alpha}
\bibliography{Biblio}

\end{document}